\documentstyle[12pt]{article}

\def\O{{\cal O}}

\newtheorem{theorem}{Theorem}[section]
\newtheorem{proposition}[theorem]{Proposition}
\newtheorem{lemma}[theorem]{Lemma}
\newtheorem{definition}[theorem]{Definition}

\def\rem{\refstepcounter{theorem}\paragraph{Remark \thetheorem}}

\def\proof{\paragraph{Proof}}

\makeatletter \def\l@section{\@dottedtocline{1}{0em}{1.2em}} \makeatother

\title{Quasi-parabolic Siegel Formula}
\author{Nitin Nitsure}
\begin{document}
\date{Corrected version, 6 December 1996}
\maketitle
School of Mathematics, Tata Institute of Fundamental Research,
Homi Bhabha Road, Mumbai 400 005, India. e-mail:
nitsure@math.tifr.res.in

\begin{abstract} The result of Siegel that the Tamagawa number
of $SL_r$ over a function field is $1$ has an expression purely
in terms of vector bundles on a curve, which is known as the
Siegel formula. We prove an analogous formula for vector bundles
with quasi-parabolic structures. This formula can be used to
calculate the betti numbers of the moduli of parabolic vector
bundles using the Weil conjucture.
\end{abstract}

\section{Introduction}
The Betti numbers of the moduli of stable vector bundles on a
complex curve, in all the
cases where the rank and degree are coprime,
were first determined by Harder and Narasimhan [H-N] as an
application of the Weil conjuctures. For this, they made use of
the result of Siegel that the Tamagawa number of the special
linear group over a function field is 1. In their refinement of
the same Betti number calculation in [D-R], Desale
and Ramanan expressed the result of Siegel in purely vector
bundle terms. This result about the Tamagawa number, called the
Siegel formula, was later given a simple proof in the language
of vector bundles by Ghione and Letizia [G-L], by introducing a
notion of effective divisors of higher rank on a curve, and
counting the number of effective divisors which correspond to a
given vector bundle. This purpose of this note is to introduce
the notion of a quasi-parabolic divisor of higher rank on a curve
(Definition 3.1 below), and to prove a quasi-parabolic analogue
(Theorem 3.4 below) of the Siegel formula, which is done here by
suitable generalizing the method of [G-L]. In a note to follow,
this formula is used to calculate the Zeta function and thereby
the Betti numbers of the moduli of parabolic bundles in the case
`stable = semistable' (these Betti numbers have already been
calculated by a guage theoretic method for genus $\ge 2$ in [N]
and for genus $0$ and $1$ by Furuta and Steer in [F-S]).

{\bf Acknowledgement} I thank M. S. Narasimhan for suggesting
the problem of extending [H-N] to parabolic bundles.

\section{Divisors supported on $X-S$}

Let $X$ be an absolutely irreducible, smooth projective curve
over the finite field $k={\bf F}_q$, and let $S$ be any closed
subset of $X$ whose points are $k$-rational. Let $K$ denote the
function field of $X$, and let $K_X$ denote the constant sheaf $K$ on
$X$. Let $g$ denote the genus of $X$. Let $r$ be a positive
integer. 
Recall that (see [G-L]) a coherent subsheaf $D\subset K _X^r$ of
generic rank $r$ is called an $r$-divisor, and the $r$-divisor
is called effective (or positive) if $\O _X^r\subset D$. The
support of the divisor is by definition 
the support of the quotient $D/\O _X^r$, which is a torsion sheaf. 
The lenght $n$ of $D/\O _X^r$ is called the degree of the divisor. 
Note that $D$ is a locally free sheaf of rank $r$ and degree
$n$. 

\rem  %2.1
Let $Z_X(t)$ be the zeta function of $X$. Then as $S$ consists
of $k$-rational points, it can be seen that the zeta function
$Z_{X-S}$ of $X-S$ is given by the formula
$$Z_{X-S}(t) = (1-t)^sZ_X(t)\eqno(1)$$
where $s$ is the cardinality of $S$. 

Note that an effective $r$-divisor on $X-S$ is the same as an
effective $r$-divisor on $X$ whose support is disjoint from $S$.
The part (1) of the proposition 1 of [G-L] gives the following,
with $X-S$ in place of $X$.

\begin{proposition} %2.2
Let $b_n^{(r)}$ be the number of effective $r$-divisors of degree
$n$ on $X$ whose support is disjoint from $S$. Let
$Z_{X-S}^{(r)}(t) = \sum _{n\ge 0} b_n^{(r)}t^n$. 
Then we have 
$$Z_{X-S}^{(r)}(t) = \prod _{1\le j\le r} Z_{X-S}(q^{j-1}t)\eqno(2)$$
\end{proposition}

In order to have the analogue of the part (2) of the proposition
1 of [G-L], we need the following lemmas.

\begin{lemma} %2.3
Let $V$ be a finite dimensional vector space over $k={\bf F}_q$,
and $s$ a positive integer. For any $1\le i\le s$, 
let $\pi _i:k^s\to k$ be the linear 
projection. For any surjective linear map $\phi :V\to k^s$, let
$V_i$ be the kernel of $\pi _i\phi :V\to k$, which is a
hyperplane in $V$ as $\phi$ is surjective. 
Let $P=P(V)$, and $P_i=P(V_i)$ denote the corresponding
projective spaces. Let $N(\phi )$ denote the number of
$k$-rational points of $P - \cup _{1\le i\le s} P_i$. Then for
any other surjective $\psi : V\to k^s$, 
we have $N(\phi )=N(\psi )$. In other words, given $s$, this
number depends only on $dim(V)$. 
\end{lemma}

\proof Given any two surjective maps $\phi ,\psi :V\to k^s$, there
exists an $\eta \in GL(V)$ such that $\phi \eta = \psi$. From
this, the result follows.

\begin{lemma} %2.4
Let $n$ be a positive integer, such that $n>2g-2+s$ 
where $g$ is the genus of $X$ and $s$ is the cardinality of $s$.
Let $b_n$ is the total number
of effective $1$-divisors of degree $n$ supported on $X-S$. 
Then for any line bundle $L$ on $X$ of degree $n$, the number of effective
$1$-divisors supported on $X-S$ which define $L$ is $b_n/P_X(1)$, where
$P_X(1)$ is the number of isomorphism classes of line bundles of
any fixed degree on $X$. 
\end{lemma}

(Here, $P_X(t)$ is the polynomial $(1-t)(1-qt)Z_X(t)$.)

\proof Let $L$ be any line bundle on $X$ of degree $n$, where
$n>2g-2+s$. Then $H^1(X,L(-S))=0$, so the natural map 
$\phi :H^0(X,L)\to H^0(X,L|S)$ is surjective. Let $V=H^0(X,L)$.
Then $dim(V)=n+1-g$. Choose a basis for each fiber $L_P$ for
$P\in S$. This gives an identification of $H^0(X,L|S)$ with
$k^s$. Now it follows that the number $N(\phi )$ defined in the
preceeding lemma depends only on $n$, and is independent of the
choice of $L$ as long as it has degree $n$. But $N(\phi )$ is
precisely the number of effective $1$-divisors supported on $X-S$,
which define the line bundle $L$ on $X$. 

Using the above lemma, the following proposition follows, by an
argument similar to the proof of part (2) of proposition 1 in
[G-L]. The proof in [G-L] expresses the number
of $r$-divisors in terms of the number of $1$-divisors, and 
the above lemma tells us the number of $1$-divisors with support
in $X-S$ corresponding to a given line bundle on $X$.

\begin{proposition} %2.5
For $L$ a line bundle of degree $n$, let $b_n^{(r,L)}$ be the 
number of effective $r$-divisors on $X$ supported on $X-S$,
having determinant isomorphic to $L$. Then provided that
$n>2g-2+s$, we have 
$$b_n^{(r,L)} = b_n^{(r)}/P_X(1)\eqno(3)$$
\end{proposition}

\begin{proposition} %2.6
$$\lim _{n\to\infty} {b_n^{(r)}\over{q^{rn}}} = 
P_X(1){(q-1)^{s-1}\over{q^{g-1+s}}}
Z_{X-S}(q^{-2})\cdots Z_{X-S}(q^{-r})\eqno(4)$$
\end{proposition}

\proof The above statement is the analogue of proposition 2 of
[G-L], with the following changes. Instead of all $r$-divisors
on $X$ in [G-L], we consider only those which are supported over
$X-S$, and instead of $Z_X(t)$, we use $Z_{X-S}(t)$. As
$Z_{X-S}(t) = (1-t)^sZ_X(t)$, the property of $Z_X(t)$ that it
has a simple pole at $t=q^{-1}$ and is regular at $1/q^j$ for
$j\ge 2$ is shared by $Z_{X-S}(t)$. 
Hence the proof in [G-L] works also in our case, proving the
proposition. 

\rem %2.7
There is a minor misprint in the equation labeled (1) in [G-L]
(page 149); the factor $q^{g-1}$ should be read as $q^{1-g}$. 

Let $L$ be any given line bundle on $X$. 
Choose any closed point $P\in X-S$, and let $l$ denote its
degree. For any $\O _X$ module $E$, set $E(m)=E\otimes
\O_X(mP)$. If a vector bundle $E$ of rank $r$ degree $n$ has
determinant $L$, then $E(m)$ has determinant $L(rm)$, degree
$n+rml$ and Euler characteristic $\chi(m)=n+rml+r(1-g)$. 

The equations (3) and (4) above imply the following.

$$\lim _{m\to\infty} 
{b_{n+rml}^{(r,L(rm))}\over{q^{r\chi(m)}}} = 
(q-1)^{s-1} 
q^{(r^2-1)(g-1)-s} 
Z_{X-S}(q^{-2})\cdots Z_{X-S}(q^{-r})\eqno(5)$$

\section{Quasi-parabolic divisors}
For basic facts about parabolic bundles, see [S] and [M-S]. We
now introduce the notion of a quasi-parabolic effective divisor
of rank $r$. Let $S\subset X$ be a finite subset consisting of
$k$-rational points. For each $P_i\in S$, let there be given
positive integers $p_i$ and $r_{i,1},\ldots ,r_{i,p_i}$ with
$r_{i,1}+\ldots +r_{i,p_i} =r$. This will be called, as usual,
the quasi-parabolic data. Recall that a quasi-parabolic
structure on a vector bundle $E$ of rank $r$ on $X$ by
definition consists of flags 
$E_{P_i}=F_{i,1}\supset\ldots\supset F_{i,p_i}\supset F_{i,p_i+1}=0$ 
of vector subspaces in the fibers over the points of $S$ such that
$dim(F_{i,j}/F_{i,j+1})=r_{i,j}$ for each $j$ from $1$ to $p_i$.  

\begin{definition}\rm %3.1
Let $X$, $S$, and the numerical data $(r_{i,j})$ be as above.
A positive quasi-parabolic divisor $(F,D)$ on $X$ consists of 
(i) a quasi-parabolic structure $F$ on the trivial bundle $\O
_X^r$, consisting of flags $F_i$ in $k^r$ at points $P_i\in S$
of the given numerical type $(r_{i,j})$, 
together with 
(ii) an effective $r$-divisor $D$ on $X$, supported on $X-S$.  
\end{definition}

Note that if $(F,D)$ is a quasi-parabolic $r$-divisor, then the
rank $r$ vector bundle $D$ has a parabolic structure given by
$F$. We denote by $P^{(r)}_E$ the set of all effective parabolic
$r$-divisors whose associated parabolic bundle is isomorphic to
a given parabolic bundle $E$.
For any vector bundle $E$ of rank $r$, let
$Hom^S_{inj}(\O _X^r,E)$ denote the set of all injective
sheaf homomorphisms $\O _X^r \to E$ which are injective when
restricted to $S$. For any quasi-parabolic bundle $E$, the group
of all quasi-parabolic automorphisms of $E$ will be denoted by
$ParAut(E)$. Then $ParAut(E)$ acts on $Hom_{inj}^S(\O _X^r,E)$
by composition. This action is free, and $P^{(r)}_E$ has a 
canonical bijection with the quotient set $Hom_{inj}^S(\O
_X^r,E)/ParAut(E)$. Hence the cardinality of $P^{(r)}_E$ is
given by
$$|P^{(r)}_E|=
{{|Hom_{inj}^S(\O _X^r,E)|}\over |{ParAut(E)|}}\eqno(6)$$

For $1\le i\le s$, let ${\rm Flag}_i$ be the variety of flags in
$k^r$ of the numerical type $(r_{i,1},\ldots ,r_{i,p_i})$. 
Let ${\rm Flag}_S = \prod _{1\le i\le s}{\rm Flag}_i$. Let
$f(q,r_{i,j})$ denote the number of $k$-rational points of ${\rm
Flag}_S$. If $a_n^{(r,L)}$ denotes the number
of quasi-parabolic divisors of flag data $(r_{i,j})$ with
degree $n$, rank $r$ and determinant $L$, then we have
$$a_n^{(r,L)} = f(q,r_{i,j})b_n^{(r,L)}\eqno(7)$$

Now let $J(r,L)$ denote the set of all isomorphism classes of
quasi-parabolic vector bundles of rank $r$, degree $n$, determinant
$L$ having the given quasi-parabolic data $(r_{i,j})$ over $S$. Hence
the equation (6) above implies the following.
$$a^{(r,L)}_n = \sum _{E\in J(r,L)}
{{|Hom_{inj}^S(\O _X^r,E)|}\over |{ParAut(E)|}}\eqno(8)$$
For any integer $m$, the map from $J(r,L) \to J(r, L(rm)$ which sends 
$E$ to $E(m) =E\otimes O_X(mP)$ is a bijection which preserves
$|ParAut|$. Hence for each $m$, we have
$$a^{(r,L(rm))}_{n+rml} = \sum _{E\in J(r,L)}
{{|Hom_{inj}^S(\O _X^r,E(m))|}\over |{ParAut(E)|}}\eqno(9) $$

\begin{lemma}
With the above notations,
$$\lim_{m\to \infty}{{|Hom_{inj}^S(\O _X^r,E(m))|}\over {q^{r\chi (E(m))}}} 
= {(q^r-1)^s(q^r-q)^s\cdots(q^r-q^{r-1})^s\over q^{r^2s}}\eqno(10)$$
If $S$ is non-empty, the limit is already attained for all 
large enough $m$ (where `large enough' depends on $E$).
\end{lemma}

\proof If $S$ is empty, the above lemma reduces to 
lemma 3 in [G-L]. If $S$ is nonempty, then any morphism of locally free sheaves
on $X$ which is injective when restricted to $S$ is injective. 
Let $m$ be large enough, so that $E(m)$ is generated by
global sections, $H^1(X,E(m))=0$, and $h^0(X,E(m))=\chi(E(m)) \ge rs$. 
Then $H^0(X,E(m))$ has a basis consisting of 
sections $\sigma_{i,P_j}$, $\tau_{\ell}$ for $i=1,\ldots, r$, $j=
1,\ldots,s$, and $\ell=1,\ldots, \chi(E(m))-rs$, such that

(1) the sections $\tau_{\ell}$ are zero on $S$,

(2) the sections $\sigma_{i,P_j}$ are zero at all other points of $S$
except $P_j$ (and hence $\sigma_{i,P_j}$ restrict at $P_j$ to a basis
of the fiber of $E(m)$ at $P_j$.

Any element of $Hom_{\O_X}(\O_X^r,E(m))
= Hom_{{\bf F}_q}({\bf F_q}^r, H^0(X,E(m)))$ is given in terms of this
basis by a $r\times q^{\chi(E(m))}$ matrix $A$. The condition that
this lies in

$$Hom_{inj}^S(\O _X^r,E(m)) \subset Hom(\O _X^r, E(m))$$
 
is the condition that each of the $s$ disjoint $r\times r$-minors,
corresponding to the part $\sigma_{1,P_j},\ldots, \sigma_{r,P_j}$ of
the basis, has nonzero determinant. This contributes the factor 
$${|GL_r({\bf F}_q)|\over |M_r({\bf F}_q)|}=
{(q^r-1)(q^r-q)\cdots(q^r-q^{r-1})\over q^{r^2}} 
$$ 
for each $P_j$, which
proves the lemma.

\begin{lemma} % 3.3
The following sum and limit can be interchanged to give
$$\sum _{E\in J(r,L)} \lim _{m\to\infty} 
{{|Hom_{inj}^S(\O _X^r,E(m))|}\over {q^{r\chi (E(m))}|ParAut(E)|}}
= \lim _{m\to\infty} \sum _{E\in J(r,L)}
{{|Hom_{inj}^S(\O _X^r,E(m))|}\over {q^{r\chi (E(m))}
|ParAut(E)|}}
$$ 
\end{lemma}

This lemma has a proof entirely analogous to the corresponding
statement in [G-L], so we omit the details.

By equation (10), the left hand side in the above
lemma equals 
$${(q^r-1)^s(q^r-q)^s\cdots(q^r-q^{r-1})^s\over q^{r^2s}}
\sum _{E\in J(r,L)} {1\over{|ParAut(E)|}}
$$

On the other hand, by (9), the right hand side is  
$\lim _{m\to\infty} a^{(r,L(rm))}_{n+rml}/q^{r\chi (m)}$. By
equations (5) and (7), this limit has the following value.  
$$f(q,r_{i,j})(q-1)^{s-1} q^{(r^2-1)(g-1)-s} 
Z_{X-S}(q^{-2})\cdots Z_{X-S}(q^{-r})$$
By putting $Z_{X-S}(t) = (1-t)^sZ_X(t)$ in the above, and cancelling
common factors from both sides, we get the following.

\begin{theorem}%3.4 
{\rm (Quasi-parabolic Siegel formula)}
$$\sum _{E\in J(r,L)} {1\over{|ParAut(E)|}}=
f(q,r_{i,j}) {q^{(r^2-1)(g-1)}\over q-1 } 
Z_X(q^{-2})\cdots Z_X(q^{-r}) $$
\end{theorem}

\rem %3.5
If $S$ is empty or more generally if the quasi-parpbolic structure at 
each point of $S$ is trivial (that
is, each flag consists only of the zero subspace and the whole space), 
then on one hand $ParAut(E)=Aut(E)$, and on the other hand each flag 
variety is a point, and so $f(q,r_{i,j})=1$. Hence in this situation 
the above formula reduces to the original Siegel formula
$$\sum _{E\in J(r,L)} {1\over{|Aut(E)|}}=
{q^{(r^2-1)(g-1)}\over q-1 } Z_X(q^{-2})\cdots Z_X(q^{-r}) $$

\section*{References} 
[D-R] Desale, U. V. and Ramanan, S. : Poincar\'e Polynomials of
the Variety of Stable Bundles, {\sl Math. Annln.} {\bf 216}
(1975), 233-244.

[F-S] Furuta, M. and Steer, B. : Siefert-fibered homology
3-spheres and Yang-Mills equations on Riemann surfaces with
marked points, {\sl Adv. Math.} {\bf 96} (1992) 38-102. 

[G-L] Ghione, F. and Letizia, M. : Effective divisors of higher
rank on a curve and the Siegel formula, {\sl Composito Math.}
{\bf 83} (1992), 147-159.

[H-N] Harder, G. and Narasimhan, M. S. : On the Cohomology
Groups of Moduli Spaces of Vector Bundles over Curves, {\sl
Math. Annln.} {\bf 212} (1975), 215-248.

[M-S] Mehta, V. B. and Seshadri, C. S. : Moduli of vector bundles
on curves with parabolic structures, {\sl Math. Annln.} {\bf 248}
(1980) 205-239.

[N] Nitsure, N. : Cohomology of the moduli of parabolic vector
bundles, {\sl Proc. Indian Acad. Sci. (Math. Sci.)} {\bf 95}
(1986) 61-77.

[S] Seshadri, C. S. : Fibres vectoriels sur les courbes
algebriques, {\sl Asterisque} {\bf 96} (1982).

%\bigskip
%Address:

%School of Mathematics, Tata Institute of Fundamental Research,
%Homi Bhabha Road, Bombay 400 005, India. e-mail:
%nitsure@tifrvax.tifr.res.in

\end{document}